
%
\input amstex
\documentstyle{amsppt}
\magnification\magstep1
\nologo
\NoBlackBoxes
\NoRunningHeads
\def\wwY#1{\widehat{Y}^{#1}}
\def\oo{\varnothing}
\def\barCP#1{\overline{\C\roman P}^{#1}}
\def\conj{\mathop{\roman{conj}}\nolimits}
\def\e{\varepsilon}
\def\z2{\Z/2}

\def\ras{$(X,\conj)$ }
\def\xr{X_\R }
\def\o{\Cal O}


\def\rank{\mathop{rank}}
 
 \def\ttx#1{\widehat{X}_0^{(#1)}}
\def\on{\overline}
\def\f{\overline{\varphi}}
\def\W{\widehat V}
\pagewidth{28pc}
\pageheight{43pc}
\topskip=\normalbaselineskip
\let\ge\geqslant

\def\C{{\Bbb C}}
\def\R{{\Bbb R}}
\def\Z{{\Bbb Z}}

\def\Cp#1{\C\roman P^{#1}}
\def\CP#1{\C\roman P^{#1}}

\let\tm\proclaim
\let\endtm\endproclaim

\let\rk=\remark
\let\endrk=\endremark

\topmatter
\title
Quotients
by complex conjugation
for Real complete intersection surfaces
\endtitle
\author
S\. M\. Finashin
\endauthor
\address
Middle East Technical University,
Ankara, 06531, Turkey
\newline\indent
St\.-Petersburg Electrotechnical University,
199376, Russia
\endaddress{}
\email serge\,\@\,rorqual.cc.metu.edu.tr
\endemail
\endkeywords
\endsubjclass
\abstract
Quotients $Y=X/\conj$ by the complex conjugation
$\conj\: X\to X$
 for complex surfaces $X$ defined over $\R$
tend
to be completely decomposable when they are simply connected,
i.e., split into connected sums
$\#_n \Cp2\#_m\barCP2$ if $w_2(Y)\ne0$,
or into $\#_n(S^2\!\times\! S^2)$ if  $w_2(Y)=0$.
The author proves this property for
 complete intersections  which are constructed
by method of a small perturbation.
\endabstract
\endtopmatter
%
%
\document
\heading
\S 1. Introduction
\endheading
We mean by a Real variety (Real curve, Real surface etc.)
a pair $(X,\conj)$, where
$X$ is a complex variety and $\conj\:X\to X$
 an anti-holomorphic involution called
{\it the real structure} or {\it the complex conjugation}.
Given an algebraic variety over $\R$
we consider the set of its complex points with the natural complex
conjugation (the Galois transformation) as the corresponding
Real variety.
The fixed point set of $\conj$ will be denoted by $\xr$ and
called {\it the real part of $X$}.
We put $Y=X/\!\conj$ and confuse in what follows
$\xr$ with its image
$q(\xr)$ under the quotient map $q\:X\to Y$.

If \ras is a nonsingular Real curve then
$Y$ is a compact surface with the boundary
$\xr$ and its topological type depends, obviously, only on
the genus of $X$, the number of components in $\xr$ and 
orientability of $Y$; the latter depends on
vanishing of the fundamental
class $[\xr]\in H_1(X;\z2)$. It is also trivial enough that any compact
surface with nonempty boundary can appear as $Y$ for some Real
curve.

The subject of the author's interest is topology of $Y$ in the case
of nonsingular Real surfaces. It is not difficult to see that
$Y$ in this case is a closed 4-manifold with the map
$q$ being a $2$-fold covering branched along $\xr$. Moreover,
$Y$ inherits from $X$ an orientation and
a smooth structure 
making $q$ smooth and orientation preserving.
The natural question is to describe the diffeomorphism types
of 4-manifolds which can arise as the quotients
$Y$ for Real  surfaces $X$.
Another interesting question
is \  if the topological types of $X$ and $\xr$ together with
some information about
the fundamental class $[\xr]\in H_2(X;\z2)$
determine the topology of
$Y$ like in the case of curves.

It turns out that
in all the known examples
if $Y$ is simply  connected then it splits
into a connected sum of copies of
$\CP2$, $\barCP2$ or $S^2\times S^2$.
Let call this property CDQ-property
(complete decomposability for quotients)
and call the Real surfaces satisfying it
 CDQ-surfaces.
Note that $Y$ is simply  connected, for example,
if $X$ is simply  connected
and   $\xr\ne\oo$. 

CDQ-property is known
for all Real rational  \cite{3} and all Real $K3$
surfaces \cite{2} except the cases when $X_{\R}=\oo$.
Further, CDQ-property is proved
for a plenty of double planes,
for a family of elliptic surfaces
and in fairly general setting for doubles
of CDQ-surfaces branched along certain ``double'' curves
(for details and more history see \cite{2}).

In the present paper we show that CDQ-surfaces
can be found among complete intersections in $\Cp{n+2}$ of arbitrary
 multi-degrees $(d_1,\dots,d_n)$.
More precisely, it appears that the straightforward
 construction of Real complete intersections
by a small perturbation method
produces CDQ-surfaces.

This gives some arguments for  the following
relatively moderate conjecture:
{\it any deformation type of simply  connected complex algebraic surfaces
contains a CDQ-surface.}



\heading
\S2. Main results
\endheading

Let $(V,\conj)$ be a Real variety.
A holomorphic linear bundle  $p\: L\to V$
 will be called {\it a Real bundle} if it
 is supplied with an anti-linear involution
$\conj_{L}\: L \to L$, which covers $\conj$:

$$
\CD
L     @>\conj_{L}>>     L\\
@VpVV              @VpVV\\
V      @>\conj>>       V\\
\endCD
$$

A section $f\:V\to L$ is called real if
$\conj_{L}\circ f=f\circ\conj$. The zero divisor 
of $f$ is, clearly, a Real  subvariety of $V$.

Assume now that $(V, \conj)$ is a Real nonsingular and connected 3--fold,
$L_i\to V$, $i=1,2$, are real linear bundles and
$f_i\:V\to L_i$ are real sections which zero divisors $X^{(i)}_0$ are
nonsingular and intersecting  transversally. Assume further that
$f\:V\to L$ is a real section of $L=L_1\otimes L_2$ which
zero divisor $X$ intersects transversally the surfaces
 $X_0^{(i)}$, $i=1,2$, and
the curve $A=X_0^{(1)}\cap X_0^{(2)}$.
Consider the section $f_\e\:V\to L_0$, $f_\e=f_1\otimes f_2+\e f$,
 $\e\in \R$,
and denote by $X_\e$ its
zero divisor, which is nonsingular
for a sufficiently small $\e> 0$,
as it can be easily seen.

\tm{Theorem 1}
If $X_0^{(1)}$ and $X_0^{(2)}$ are CDQ-surfaces,
$A$ is connected, $A_\R\ne\varnothing$
 and $\e>0$ is small enough then
 $X_\e$ is CDQ-surface as well.
\endtm

We prove this theorem in \S 3 and discuss in the rest of
this section some of its
implications.

Let $(V,\conj)$ be as above.
We call a real linear bundle
 $L$ CDQ{\it-bundle} if it is  very ample and admits
a real section with a nonsingular
CDQ zero divisor.

\tm{Lemma 2}
If $L$ is a CDQ-bundle then its multiples
$L^{\otimes d}$, $d\ge 1$, are
CDQ-bundles as well.
\endtm

\demo{Proof}
Let $X_0^{(1)}$ be a CDQ-divisor of $L$.
We prove by induction on $d$ that there exists
a CDQ-divisor, $X_0^{(2)}$,
of $L^{\otimes d}$ which intersects
$X_0^{(1)}$
transversally along a curve having nonempty real part.
This claim is trivial for $d=1$, since
we can perturb $X_0^{(1)}$ so that the result will
intersect  $X_0^{(1)}$ transversally and contain 
a given real point of it.

Suppose that $X_0^{(2)}$ satisfies the induction assumption.
By Lefschetz Theorem $A$ is connected.
A generic real section of  $L^{\otimes(d+1)}$ has zero divisor $X$
transversal to $X_0^{(i)}$ and to $X_0^{(1)}\cap X_0^{(2)}$,
hence, we can apply Theorem 1 and get a CDQ-divisor $X_\e$
by a perturbation of $X_0^{(1)}\cup X_0^{(2)}$ via $X$.
We can also choose $X$ containing a real point of $X_0^{(1)}$,
since $L^{\otimes(d+1)}$ is very ample.
Then, for a sufficiently small $\e>0$,  $X_\e$
 intersects $X_0^{(1)}$ transversally and
 $X_{\e}\cap X_0^{(1)}=X\cap X_0^{(1)}$ has nonempty real part.
\qed
\enddemo

\tm{Theorem 3}
For arbitrary integers $n, d_1,\dots,d_n\ge 1$
there exists a CDQ-surface $X\subset\Cp{n+2}$ which is
a complete intersection of multi-degree $(d_1,\dots,d_n)$.
\endtm
\demo{Proof}
Induction on $n$. $\Cp2/\conj$ is
diffeomorphic to $S^4$ (see, e.g., \cite{1}),
therefore,
 $\o(1)_{\Cp3}$  is a CDQ-bundle.
By Lemma 2, $\o(d)_{\Cp3}$, $\forall d\ge 1$, is a CDQ-bundle as well.
Assume now that we are given
 a complete intersection of Real hypersurfaces,
$X=H_1\cap\dots\cap H_n\subset\Cp{n+2}$,
of multi-degree $(d_1,\dots,d_n)$ and that $X$ is CDQ-surface.
Choose  hypersurfaces
$H'_i\subset\Cp{n+3}$, $i=1,\dots, n$,
so that
$H_i=H'_i\cap\Cp{n+2}$,
and the intersection
$V=H'_1\cap\dots\cap H'_n$ is transversal.
Then the bundle $L\to V$ induced from $\o(1)_{\Cp{n+3}}$
is CDQ-bundle, since $X$ is its zero divisor. By Lemma~2, $L^{\otimes d}$
is also CDQ-bundle, hence, there exists
a CDQ complete intersection of
multi-degree $(d_1,\dots,d_n,d)$.
\qed
\enddemo

\rk{Remark}
The method used for Theorem 3 can be applied similarly
to complete intersections
in weighted projective spaces or in products of projective spaces
and yields also
CDQ-surfaces of arbitrary multi-degree.
\endrk

\heading
\S3 Proof of Theorem 1
\endheading

Denote by $N^{(i)}$  a $\conj$-invariant
tubular neighborhood of $A$ in $X_0^{(i)}$ and
put $\overline N^{(i)}=N^{(i)}/\conj$, $B=A/\conj$, $Y^{(i)}=X_0^{(i)}/\conj$
and $Y_\e=X_\e/\conj$.
It can be easily seen that $\overline N^{(i)}$ is a regular
neighborhood of $B$ in $Y^{(i)}$.
Let $2k$ denote the number of
 imaginary points in
$A\cap X$.

\tm{Proposition 4}
There exists a diffeomorphism
$\f\: \partial\on N^{(1)}\to\on N^{(2)}$,
such that
$Y_\e\cong((Y^{(1)}-\on N^{(1)})\cup_{\f}(Y^{(2)}-\on N^{(2)}))\,
\#k\,\barCP2$.
\endtm

Let see first how Theorem 1 follows from the above proposition.

Since $A$ is connected and has nonempty real part,
 $B$ is a compact connected surface with a nonempty
boundary, hence, $\overline N^{(i)}$  are
 handlebodies with one $0$-handle and several $1$-handles
embedded into $Y^{(i)}$.
It is well known that if we glue
a pair of simply  connected 4-manifolds, $Y^{(i)}, i=1,2$,
along the boundary of such handlebodies the result
is diffeomorphic to $Y^{(1)}\#Y^{(2)}\#g\,Z$, where
$g=\rank(H_1(B))$ and $Z=S^2\times S^2$ or $\Cp2\# \barCP2$
(see, e.g., \cite{4}).
This implies complete decomposability
of $Y_\e$ if $Y^{(i)}$ are completely decomposable.

The proof of Proposition 4 follows closely the idea of \cite{4}.
By blowing up $\W\to V$ we lift the pencil $X_t=(1-t)X_0+tX$
to a real fibering $\W\to\Cp1$. Then we
apply the deformation theorem of \cite{4} in the equivariant version.
Specifically,
assume that $\e\in\R$  and $\e>0$ is sufficiently small.
Then $X_\e$ intersects
$X_0^{(i)}$,  $i=1,2$ transversally along the curve
 $C_i=X_0^{(i)}\cap X$.
Consider first
the blow-up, $\widetilde V\to V$, along $C_1$
and denote by $\widetilde C_i$, $\widetilde X_0^{(i)}$,
and $\widetilde X_t$ the proper images of
$C_i$, $X_0^{(i)}$ and $X_t$.
The pencil $\widetilde X_t$  has the base-curve
$\widetilde C_2$, therefore, the next
blow-up $\W\to \widetilde V$ along  $\widetilde C_2$
gives a fibering over $\Cp1$ with fibers
$\widehat{X}_t$
(here and below we mark by a hat
the proper image in $\W$).

The projections $\widehat{X}_\e\to X_\e$,
$\widehat{X}_0^{(1)}\to X_0^{(1)}$
are biregular, as well as
$\widehat{X}^{(2)}_0\to\widetilde X^{(2)}_0$,
whereas
$\widetilde X^{(2)}_0\to X^{(2)}_0$
is the blow-up at $C_1\cap X_0^{(2)}=A\cap X$.

The real structure on $V$ can be obviously lifted to
the real structure,  $\conj_{\W}\:\W\to\W$,
 and we have $\wwY{(1)}\cong\ Y^{(1)}$,
$\wwY{}_\e\cong\ Y_\e$ and
$\wwY {(2)}\cong\ Y^{(2)}\#k(\barCP2)$,
where $\wwY {(1)}$,
$\wwY{}_\e$, $\wwY {(2)}$
denote the quotients by $\conj_{\W}$ of
$\ttx1$, $\widehat{X}_{\e}$ and $\ttx2$.
The latter diffeomorphism follows from that
blow-ups at real points do not change the diffeomorphism type
of the quotient, since $\Cp2/\conj\cong S^4$,
whereas a pair of blow-ups at conjugated imaginary points
descends to a blow-up in the quotient.
Restrictions
give diffeomorphisms
between $N^{(i)}$ and regular neighborhoods of
$\widehat A/\conj_{\W}$ in $\ttx{i}/\conj_{\W}$ for
$\widehat A=\ttx1\cap\ttx2$.

To complete the proof we use
the deformation
theorem \cite{4}. Recall the statement of one of its corollaries.

{\it
Let $f\: W\to \Delta$ be a nonconstant proper holomorphic mapping
of a $3$-fold $W$ into a disc,
$\Delta\subset\C$, around zero.
Assume that $f$ has
a critical  value only at zero and the zero divisor $X_0$ of $f$
consists of two nonsingular irreducible components
$X_0^{(i)}$, $i=1,2$, of multiplicity $1$ crossing transversally
along a nonsingular irreducible curve $A$.
Suppose that $U\subset W$ is a sufficiently small
tubular neighborhood of $A$,
so that $ N^{(i)}=U\cap X^{(i)}_0$
is a tubular neighborhood of $A$ in $X^{(i)}_0$, $i=1,2$.
 Then there exists a bundle
isomorphism $\varphi\:\partial N^{(1)}\to\partial N^{(2)}$,
 reversing  orientations of fibers,
such that $X_t=f^{-1}(t)$ is diffeomorphic to
$(X_0^{(1)}-N^{(1)})\cup_{\varphi}(X_0^{(2)}-N^{(2)})$
for a non-critical value $t\in\Delta$.
}

In the equivariant version of this theorem we
assume in addition that $W$ is supplied with a real structure
$\conj_W\:W\to W$ and
$f\circ\conj_W=\conj\circ f $, where $\conj\:\Delta\to\Delta$ is the
complex conjugation on $\C$.
Then one can choose
the neighborhood $U$ to be $\conj_{W}$-invariant,
make the isomorphism $\varphi$
$\conj_W$-equivariant
and the diffeomorphism
$X_t\cong(X_0^{(1)}-N^{(1)})\cup_\varphi(X_0^{(2)}-N^{(2)})$
commute with
the involutions of complex conjugation
in $X_t$ and $X_0^{(i)}-N^{(i)}$.

To prove the latter
one should repeat the arguments
of \cite{4} with some, not essential modifications:
we need to choose
a $\conj_W$--invariant metric on $W$ and
 instead of the fibering $U\cap{}X_t \to A$
considered in \cite{4}
deal with its quotient, $U\cap{}X_t/\conj_W \to B$, and then apply
similarly
the arguments on the reduction of the structure group.

\enddocument